\documentclass[12pt]{article}
\usepackage[english]{babel}
\usepackage{latexsym}
\usepackage{amssymb,amsthm,amsmath}
\usepackage{graphicx}
\title{\bf A computer based classification of caps in $PG(4,3)$}
\author{Daniele Bartoli, Stefano Marcugini, Fernanda Pambianco}

{\theoremstyle{definition}
\newtheorem*{definition*}{Definition}
}

\newtheorem*{proposition*}{Proposition}
\newtheorem*{corollary*}{Corollary}
\newtheorem*{lemma*}{Lemma}
\date{}
\begin{document}

\maketitle \vspace*{-15mm} \noindent

 \vspace*{2cm}
\begin{abstract}
\noindent In this paper we present the complete classification of
caps in $PG(4,3)$. These results have been obtained using a computer
based exhaustive search that exploits projective equivalence.
\end{abstract}

\section{Introduction}
In the projective space $PG(r,q)$ over the Galois Field $GF(q)$, a $n$-cap is a set of $n$ points
 no 3 of which are collinear. A $n$-cap is called \emph{complete} if it is not contained in a
 $(n+1)$-cap. For  a detailed description of the most important properties of these geometric
 structures, we refer the reader to \cite{Hir}.
 In the last decades the problem of determining the spectrum of the
 sizes of complete caps has been the subject of a lot of researches.
 For a survey see \cite{18secondo}.\\
In this work we search for the classification of complete and
incomplete caps in $PG(4,3)$, using an exhaustive search algorithm. In Section \ref{alg}
the algorithm utilized is illustrated;
 in Section \ref{risultati} the complete list of non equivalent complete and incomplete
caps is presented.

\section{The searching algorithm }\label{alg}
In this section the algorithm utilized is presented. Our goal is to obtain the classification of complete and incomplete caps in $PG(4,3)$. It is not restrictive to suppose that a cap in $PG(4,3)$ contains this three points:
$$\mathcal{R}=\{(1:0:0:0:0);(0:1:0:0:0);(0:0:1:0:0)\}$$
Then we define the set \emph{Cand} of all the points lying no 2-secant of $\mathcal{R}$. We introduce in \emph{Cand} the following equivalence relationship:
\begin{displaymath}
P \sim Q \iff \mathcal{C} \cup \{P\} \cong \mathcal{C} \cup \{Q\},
\end{displaymath}

\noindent where $\cong$ means that the two sets are projectively equivalent. This relationship spreads the candidates in equivalent classes $\mathcal{C}_{1},\ldots,\mathcal{C}_{k}$.\\
The choice of the next point to add to the building cap can be made only among the representatives of the equivalent classes, in fact two caps one containing $\mathcal{C} \cup \{P\}$ and the other one $\mathcal{C} \cup \{Q\}$, with $P$ and $Q$ in $\mathcal{C}_{\overline{i}}$, are equivalent by definition of orbit.\\
Suppose now that we have construct all the caps containing $\mathcal{C} \cup \{P_{i}\}$, with $i \leq \overline{i}$. Considering the caps containing $\mathcal{C} \cup \{P_{j}\}$ with $\overline{i}<j$, all the points of the classes $\mathcal{C}_{k}$ with $k<\overline{i}$ can be avoided. In fact a cap containing $\mathcal{C} \cup \{P_{j}\} \cup \{\overline{P}_{k}\}$, with $\overline{P}_{k} \in \mathcal{C}_{k}$ and $k<\overline{i}$, is projectively equivalent to a cap containing $\mathcal{C} \cup \{P_{k}\} \cup \{P_{j}\}$, already studied.\\
When we add a new point to the cap, we can divide all the remaining candidates in equivalence classes, as above. Two points  $P$ and $Q$ are in relationship with the $j$-th class $\mathcal{C}_{j}$, i.e. $P \sim_{j} Q$, if $\mathcal{C} \cup \{P_{j}\} \cup \{P\}$ and $\mathcal{C} \cup \{P_{j}\} \cup \{Q\}$ are projectively equivalent. \\
At the $m$-th step of the extension process if the cap $\mathcal{C}\cup \{P_{i_{m}}\} \cup \ldots \cup \{P^{i_{1}\ldots i_{m-1}}_{i_{m}}\} \cup\{P\}$  is projectively equivalent to the cap $\mathcal{C}\cup \{P_{i_{m}}\} \cup \ldots \cup \{P^{i_{1}\ldots i_{m-1}}_{i_{m}}\} \cup\{Q\}$ with $P^{i_{1}\ldots i_{r}}_{s} \in \mathcal{C}^{i_{1}\ldots i_{r}}_{s} $, then $P$ and $Q$ are in relationship ($P \sim_{i_{1}\ldots i_{m}} Q$) and they belong to the same class $\mathcal{C}^{i_{1}\ldots i_{m}}_{m+1} $.\\
Iterating the process we can build a tree similar to the following:
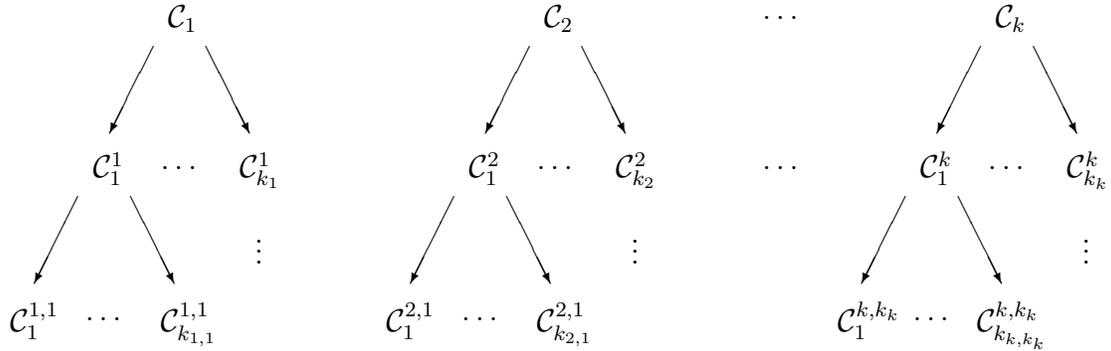
\begin{center}

\setlength{\unitlength}{1cm}
\begin{picture}(15,6)(0,0)


\put(2.35,4.9){$\mathcal{C}_{1}$}

\put(2.15,4.65){\vector(-1,-2){0.58}}
\put(2.85,4.65){\vector(1,-2){0.58}}

\put(1.35,2.9){$\mathcal{C}_{1}^{1}$}
\put(2.25,3){\ldots}

\put(3.3,2.9){$\mathcal{C}_{k_{1}}^{1}$}

\put(1.15,2.65){\vector(-1,-2){0.58}}
\put(1.85,2.65){\vector(1,-2){0.58}}

\put(3.5,1.75){\vdots}

\put(0.25,0.85){$\mathcal{C}_{1}^{1,1}$}
\put(1.25,1){\ldots}

\put(2.25,0.85){$\mathcal{C}_{k_{1,1}}^{1,1}$}

\put(10.25,5){\ldots}
\put(10.25,3){\ldots}


\put(7.35,4.9){$\mathcal{C}_{2}$}

\put(7.15,4.65){\vector(-1,-2){0.58}}
\put(7.85,4.65){\vector(1,-2){0.58}}

\put(6.35,2.9){$\mathcal{C}_{1}^{2}$}
\put(7.25,3){\ldots}

\put(8.30,2.9){$\mathcal{C}_{k_{2}}^{2}$}

\put(6.15,2.65){\vector(-1,-2){0.58}}
\put(6.85,2.65){\vector(1,-2){0.58}}

\put(8.5,1.75){\vdots}

\put(5.25,0.85){$\mathcal{C}_{1}^{2,1}$}
\put(6.25,1){\ldots}

\put(7.25,0.85){$\mathcal{C}_{k_{2,1}}^{2,1}$}


\put(13.35,4.9){$\mathcal{C}_{k}$}

\put(13.15,4.65){\vector(-1,-2){0.58}}
\put(13.85,4.65){\vector(1,-2){0.58}}

\put(12.35,2.9){$\mathcal{C}_{1}^{k}$}
\put(13.25,3){\ldots}

\put(14.3,2.9){$\mathcal{C}_{k_{k}}^{k}$}

\put(12.15,2.65){\vector(-1,-2){0.58}}
\put(12.85,2.65){\vector(1,-2){0.58}}

\put(14.5,1.75){\vdots}

\put(11.25,0.85){$\mathcal{C}_{1}^{k,k_{k}}$}
\put(12.25,1){\ldots}

\put(13.1,0.85){$\mathcal{C}_{k_{k,k_{k}}}^{k,k_{k}}$}

\end{picture}
\end{center}
The tree is important to restrict the number of candidates in the extension process. Suppose that we have generated a $n$-cap containing the cap $\mathcal{C}\cup \{P_{i_{1}}\} \cup \{P_{i_{2}}^{i_{1}}\} \cup \ldots \cup \{P^{i_{1}\ldots i_{m-1}}_{i_{m}}\} \cup\{P\}$, after having generated $n$-caps containing $\mathcal{C}\cup \{P_{j}\}$ with $j<i_{1}$, $\mathcal{C}\cup \{P_{i_{1}}\}\cup \{P^{i_{1}}_{j}\}$ with $j<i_{2}$,\ldots, $\mathcal{C}\cup \{P_{i_{1}}\}\cup \{P^{i_{1}}_{i_{2}}\}\cup \ldots \cup \{P^{i_{1}\ldots i_{m-1}}_{j}\}$ with $j<i_{m}$, with $P^{i_{1}\ldots i_{r}}_{s} \in \mathcal{C}^{i_{1}\ldots i_{r}}_{s}$. Then the points belonging to $\mathcal{C}_{1} \cup\ldots \cup \mathcal{C}_{i_{1}-1} \cup\mathcal{C}^{i_{1}}_{1} \cup\ldots \cup \mathcal{C}^{i_{1}}_{i_{2}-1} \cup \ldots \cup \mathcal{C}^{i_{1}\ldots i_{m-1}}_{1} \cup\ldots \cup \mathcal{C}^{i_{1}\ldots i_{m-1}}_{i_{m}-1}$ can be avoided, because a cap containing one of them is equivalent to one already found. For example a $n$-cap containing $\mathcal{C}\cup \{P_{i_{1}}\} \cup \ldots \cup \{P^{i_{1}\ldots i_{m-1}}_{i_{m}}\} \cup\{P\}\cup \{Q\}$ with $Q \in \mathcal{C}_{h}$ for some $h<i_{1}$ is equivalent to a $n$-cap containing $\mathcal{C} \cup \{P_{h}\}$, which is already found. \\

\section{Results}\label{risultati}
In this Section all non equivalent caps, complete and incomplete, in
$PG(5,2)$ are presented.
\subsection{Non-equivalent caps $\mathcal{K}$ in $PG(4,3)$}
This table shows the number and the type of the non equivalent
examples of all the caps.
\begin{table}[h]
\caption{Number and type of non equivalent examples}
\begin{center}
\begin{tabular}{|c|c|c||c|c|c|}
\hline
$|\mathcal{K}|$&\# COMPLETE&\# INCOMPLETE&$|\mathcal{K}|$&\# COMPLETE&\# INCOMPLETE\\
&CAPS&CAPS&&CAPS&CAPS\\
\hline
3&0&1&12&0&844\\
4&0&2&13&0&1532\\
5&0&3&14&0&2020\\
6&0&5&15&0&1778\\
7&0&8&16&48&971\\
8&0&19&17&17&320\\
9&0&46&18&32&58\\
10&0&137&19&4&16\\
11&1&355&20&9&0\\
\hline
\end{tabular}
\end{center}
\end{table}

In the following Subsections the non equivalent examples of complete
caps of sizes $11, 19$ and $20$  are presented. In particular, for
each cap we compute the stabilizer (see \cite{grouptables}) and the
weight distribution, using the software package MAGMA.

\section{11-complete cap}
\begin{center}
\begin{tabular}{ccccccccccc}
 1& 1& 0& 1& 0& 0& 0& 1& 1& 1& 0\\
 0& 0& 1& 1& 0& 0& 0& 1& 0& 0& 1\\
 0& 1& 0& 0& 1& 0& 0& 0& 1& 0& 1\\
 0& 0& 0& 0& 0& 1& 0& 0& 1& 1& 1\\
 0& 0& 0& 0& 0& 0& 1& 1& 0& 1& 1\\
\end{tabular}
\end{center}
The stabilizer is a non abelian group of order $7920$.\\
Its weight distribution is the following:
$$[<6,132>,<9,110>]$$

\section{19-complete caps}
\begin{center}
CAP 1
\end{center}
\begin{center}
\begin{tabular}{ccccccccccccccccccc}
 1& 0& 0& 0& 1& 1& 0& 0& 0& 0& 1& 1& 1& 1& 0& 1& 1& 1& 1\\
 0& 1& 0& 0& 0& 1& 0& 1& 1& 0& 1& 0& 0& 0& 1& 0& 1& 1& 0\\
 0& 0& 1& 0& 0& 0& 0& 1& 1& 1& 0& 0& 1& 1& 0& 1& 0& 0& 1\\
 0& 0& 0& 1& 1& 1& 0& 0& 1& 1& 0& 1& 0& 1& 1& 0& 1& 0& 0\\
 0& 0& 0& 0& 0& 0& 1& 1& 0& 1& 1& 1& 0& 0& 1& 1& 0& 1& 1\\
\end{tabular}
\end{center}
The stabilizer is $\mathbb{Z}_{2} \times \mathbb{Z}_{2}$.\\
Its weight distribution is the following:
$$[<10,8>,<11,72>,<12,56>,<13,38>,$$
$$<14,24>,<15,20>,<16,8>,<17,12>,<18,4>]$$

\begin{center}
CAP 2
\end{center}
\begin{center}
\begin{tabular}{ccccccccccccccccccc}
1& 0& 0& 0& 0& 1& 0& 1& 1& 1& 0& 0& 1& 1& 1& 1& 0& 1& 1\\
 0& 1& 0& 0& 0& 1& 1& 1& 1& 0& 1& 0& 1& 0& 1& 0& 1& 0& 0\\
 0& 0& 1& 0& 0& 0& 1& 0& 1& 1& 1& 1& 0& 0& 1& 1& 0& 0& 1\\
 0& 0& 0& 1& 0& 1& 0& 0& 0& 0& 1& 1& 0& 1& 0& 1& 1& 1& 0\\
 0& 0& 0& 0& 1& 0& 1& 0& 0& 0& 0& 1& 1& 1& 0& 0& 1& 1& 1\\
\end{tabular}
\end{center}
The stabilizer is a non abelian group of order $48$.\\
Its weight distribution is the following:
$$[<10,10>,<11,62>,<12,72>,<13,30>,$$
$$<14,28>,<15,2>,<16,32>,<18,6>]$$

\begin{center}
CAP 3
\end{center}
\begin{center}
\begin{tabular}{ccccccccccccccccccc}
 1& 0& 0& 0& 1& 0& 1& 1& 0& 1& 1& 0& 0& 1& 1& 0& 1& 1& 1\\
 0& 1& 0& 0& 0& 0& 0& 0& 1& 0& 1& 1& 0& 1& 0& 1& 1& 1& 0\\
 0& 0& 1& 0& 0& 0& 0& 1& 1& 0& 1& 1& 1& 1& 0& 0& 0& 0& 1\\
 0& 0& 0& 1& 1& 0& 0& 0& 0& 0& 0& 1& 1& 0& 1& 1& 0& 1& 0\\
 0& 0& 0& 0& 0& 1& 1& 0& 1& 1& 0& 0& 1& 0& 1& 1& 0& 0& 1\\
\end{tabular}
\end{center}
The stabilizer is a non abelian group of order $144$.\\
Its weight distribution is the following:
$$[<9,4>,<11,72>,<12,72>,<14,72>,<17,18>,<18,4>]$$

\begin{center}
CAP 4
\end{center}
\begin{center}
\begin{tabular}{ccccccccccccccccccc}
 1& 1& 0& 0& 0& 0& 1& 1& 0& 1& 0& 0& 1& 1& 1& 1& 0& 1& 1\\
 1& 0& 1& 0& 0& 0& 0& 1& 1& 1& 1& 0& 1& 1& 0& 0& 1& 0& 0\\
 0& 0& 0& 1& 0& 0& 0& 0& 1& 1& 1& 1& 1& 0& 0& 1& 0& 1& 1\\
 1& 0& 0& 0& 1& 0& 0& 1& 0& 0& 1& 1& 0& 0& 1& 1& 1& 1& 0\\
 0& 0& 0& 0& 0& 1& 1& 0& 1& 0& 0& 1& 0& 1& 1& 0& 1& 0& 1\\
\end{tabular}
\end{center}
The stabilizer is a non abelian group of order $48$.\\
Its weight distribution is the following:
$$[<9,2>,<10,12>,<11,48>,<12,80>,$$
$$<13,36>,<14,24>,<15,8>,<16,24>,<18,8>]$$

\section{20-complete caps}
\begin{center}
CAP 1
\end{center}
\begin{center}
\begin{tabular}{cccccccccccccccccccc}
 1& 1& 0& 1& 1& 1& 0& 0& 0& 0& 0& 1& 0& 1& 1& 1& 0& 1& 1& 1\\
 0& 0& 1& 0& 0& 1& 0& 0& 0& 1& 1& 0& 0& 1& 0& 1& 1& 1& 0& 0\\
 1& 0& 0& 0& 0& 0& 1& 0& 0& 1& 1& 0& 1& 0& 0& 1& 0& 0& 1& 1\\
 0& 0& 0& 1& 0& 0& 0& 1& 0& 0& 1& 1& 1& 0& 0& 1& 1& 0& 1& 0\\
 0& 0& 0& 1& 0& 0& 0& 0& 1& 1& 0& 0& 1& 1& 1& 1& 1& 0& 0& 0\\
\end{tabular}
\end{center}
The stabilizer is a non abelian group of order $36$.\\
Its weight distribution is the following:
$$[<11,18>,<12,96>,<13,36>,<14,36>,<15,18>,<16,18>,<18,20>]$$

\begin{center}
CAP 2
\end{center}
\begin{center}
\begin{tabular}{cccccccccccccccccccc}
 1& 1& 0& 0& 1& 0& 1& 1& 0& 1& 1& 0& 0& 1& 0& 1& 1& 1& 0& 1\\
 1& 0& 1& 0& 1& 0& 1& 0& 0& 1& 0& 1& 1& 0& 0& 1& 0& 1& 1& 0\\
 0& 0& 0& 1& 0& 0& 0& 0& 0& 1& 1& 1& 1& 0& 1& 1& 0& 0& 0& 1\\
 0& 0& 0& 0& 0& 1& 1& 0& 0& 1& 0& 0& 1& 1& 1& 0& 1& 0& 1& 0\\
 0& 0& 0& 0& 1& 0& 0& 1& 1& 1& 0& 1& 0& 0& 1& 0& 1& 0& 1& 0\\
\end{tabular}
\end{center}
The stabilizer is a non abelian group of order $36$.\\
Its weight distribution is the following:
$$[<10,6>,<11,6>,<12,96>,<13,42>,$$
$$<14,42>,<15,18>,<16,6>,<17,6>,<18,20>]$$

\begin{center}
CAP 3
\end{center}
\begin{center}
\begin{tabular}{cccccccccccccccccccc}
  1& 1& 0& 1& 1& 0& 0& 0& 0& 0& 0& 1& 1& 1& 1& 0& 1& 1& 1& 1\\
 0& 0& 1& 0& 0& 0& 0& 0& 1& 1& 0& 0& 1& 1& 1& 1& 0& 0& 0& 0\\
 1& 0& 0& 0& 0& 1& 0& 0& 1& 1& 1& 0& 0& 0& 1& 0& 1& 0& 1& 1\\
 0& 0& 0& 1& 0& 0& 1& 0& 0& 1& 1& 1& 0& 1& 1& 1& 1& 0& 0& 0\\
 0& 0& 0& 1& 0& 0& 0& 1& 1& 0& 1& 0& 1& 0& 1& 1& 0& 1& 0& 1\\
\end{tabular}
\end{center}
The stabilizer is a non abelian group of order $144$.\\
Its weight distribution is the following:
$$[<10,4>,<11,18>,<12,78>,<13,48>,$$
$$<14,36>,<15,36>,<18,20>,<19,2>]$$

\begin{center}
CAP 4
\end{center}
\begin{center}
\begin{tabular}{cccccccccccccccccccc}
 1& 1& 0& 1& 0& 0& 0& 1& 0& 1& 0& 1& 1& 1& 0& 1& 1& 1& 1& 0\\
 0& 0& 1& 0& 0& 0& 0& 1& 1& 0& 0& 1& 1& 1& 1& 0& 0& 0& 0& 1\\
 0& 1& 0& 0& 1& 0& 0& 0& 1& 1& 1& 0& 0& 0& 0& 1& 1& 0& 0& 1\\
 0& 1& 0& 1& 0& 1& 0& 1& 0& 0& 1& 1& 0& 0& 1& 0& 0& 0& 1& 1\\
 0& 0& 0& 1& 0& 0& 1& 0& 1& 1& 1& 1& 0& 1& 1& 0& 1& 1& 0& 0\\
\end{tabular}
\end{center}
The stabilizer is a non abelian group of order $2880$.\\
Its weight distribution is the following:
$$[<11,40>,<12,60>,<14,120>,<18,20>,<20,2>]$$

\begin{center}
CAP 5
\end{center}
\begin{center}
\begin{tabular}{cccccccccccccccccccc}
1& 0& 0& 0& 1& 0& 1& 1& 0& 1& 0& 0& 1& 1& 0& 1& 1& 1& 1& 1\\
 0& 1& 0& 0& 1& 0& 0& 1& 1& 1& 1& 0& 0& 1& 1& 0& 0& 0& 0& 0\\
 0& 0& 1& 0& 0& 0& 0& 0& 1& 1& 1& 1& 0& 0& 0& 0& 0& 1& 1& 1\\
 0& 0& 0& 1& 1& 0& 1& 0& 0& 0& 1& 1& 1& 0& 1& 0& 1& 0& 0& 0\\
 0& 0& 0& 0& 0& 1& 0& 0& 1& 0& 0& 1& 1& 0& 1& 0& 1& 1& 0& 1\\
\end{tabular}
\end{center}
The stabilizer is a non abelian group of order $160$.\\
Its weight distribution is the following:
$$[<12,150>,<15,72>,<18,20>]$$

\begin{center}
CAP 6
\end{center}
\begin{center}
\begin{tabular}{cccccccccccccccccccc}
1& 1& 1& 0& 1& 0& 1& 0& 0& 1& 0& 0& 1& 0& 1& 0& 1& 1& 1& 1\\
 1& 0& 0& 1& 0& 0& 1& 0& 0& 0& 1& 1& 0& 0& 0& 1& 1& 0& 0& 0\\
 1& 1& 0& 0& 0& 1& 0& 0& 0& 0& 1& 1& 0& 1& 0& 0& 0& 0& 0& 1\\
 0& 0& 0& 0& 1& 0& 0& 1& 0& 1& 0& 1& 1& 1& 0& 1& 0& 0& 0& 0\\
 0& 0& 0& 0& 1& 0& 0& 0& 1& 0& 1& 0& 0& 1& 1& 1& 0& 0& 1& 0\\
\end{tabular}
\end{center}
The stabilizer is a non abelian group of order $720$.\\
Its weight distribution is the following:
$$[<10,12>,<12,60>,<13,120>,<16,30>,<18,20>]$$

\newpage
\begin{center}
CAP 7
\end{center}
\begin{center}
\begin{tabular}{cccccccccccccccccccc}
1& 1& 0& 0& 1& 0& 1& 0& 1& 1& 0& 0& 1& 0& 1& 1& 1& 1& 0& 1\\
 1& 0& 1& 0& 1& 0& 0& 0& 1& 0& 1& 1& 0& 0& 1& 0& 1& 0& 1& 0\\
 0& 0& 0& 1& 0& 0& 0& 0& 1& 1& 1& 1& 0& 1& 1& 0& 0& 0& 0& 1\\
 0& 0& 0& 0& 0& 1& 0& 0& 1& 0& 0& 1& 1& 1& 0& 1& 0& 0& 1& 0\\
 0& 0& 0& 0& 1& 0& 1& 1& 1& 0& 1& 0& 0& 1& 0& 1& 0& 1& 1& 0\\
\end{tabular}
\end{center}
The stabilizer is the dihedral group $D_{8}$.\\
Its weight distribution is the following:
$$[<10,4>,<11,12>,<12,92>,<13,40>,$$
$$<14,42>,<15,20>,<16,10>,<18,22>]$$

\begin{center}
CAP 8
\end{center}
\begin{center}
\begin{tabular}{cccccccccccccccccccc}
1& 1& 1& 0& 1& 0& 0& 0& 1& 0& 0& 0& 1& 1& 0& 1& 1& 1& 1& 1\\
 1& 0& 0& 1& 0& 0& 0& 0& 0& 1& 1& 0& 1& 1& 1& 1& 0& 0& 0& 0\\
 1& 1& 0& 0& 0& 1& 0& 0& 0& 1& 1& 1& 0& 0& 0& 0& 1& 0& 0& 1\\
 0& 0& 0& 0& 1& 0& 1& 0& 1& 0& 1& 1& 0& 1& 1& 0& 1& 0& 0& 0\\
 0& 0& 0& 0& 1& 0& 0& 1& 0& 1& 0& 1& 1& 0& 1& 0& 0& 0& 1& 1\\
\end{tabular}
\end{center}
The stabilizer is a non abelian group of order $720$.\\
Its weight distribution is the following:
$$[<12,150>,<15,72>,<18,20>]$$

\begin{center}
CAP 9
\end{center}
\begin{center}
\begin{tabular}{cccccccccccccccccccc}
1& 1& 0& 1& 1& 1& 0& 0& 0& 0& 0& 0& 1& 1& 1& 1& 0& 1& 1& 1\\
 0& 0& 1& 0& 0& 1& 0& 0& 0& 1& 1& 0& 0& 1& 1& 1& 1& 1& 0& 0\\
 1& 0& 0& 0& 0& 0& 1& 0& 0& 1& 1& 1& 0& 0& 0& 1& 0& 0& 0& 1\\
 0& 0& 0& 1& 0& 0& 0& 1& 0& 0& 1& 1& 1& 0& 1& 1& 1& 0& 0& 0\\
 0& 0& 0& 1& 0& 0& 0& 0& 1& 1& 0& 1& 0& 1& 0& 1& 1& 0& 1& 0\\
\end{tabular}
\end{center}
The stabilizer is the dihedral group $D_{4}$.\\
Its weight distribution is the following:
$$[<10,2>,<11,16>,<12,90>,<13,44>,$$
$$<14,34>,<15,24>,<16,8>,<17,4>,<18,20>]$$

\end{document}